\documentclass{amsart}
\usepackage{amsthm,amssymb,mathtools}
\usepackage{graphicx}
\usepackage{enumitem,caption}
\usepackage{tikz}
\usetikzlibrary{decorations.pathreplacing,decorations.pathmorphing,decorations.markings,shapes}
\usetikzlibrary{shapes.geometric,shadows}
\usetikzlibrary{matrix,arrows}
\usepackage{pstricks,pstricks-add}
\psset{unit=1mm}
\usepackage{pinlabel}
\DeclareMathAlphabet{\mathpzc}{OT1}{pzc}{m}{it}
\usepackage{colortbl,multicol,multirow}
\usepackage{cancel}
\usepackage{pifont}
\usepackage[all]{xy}

\definecolor{lightgrey}{rgb}{.804,.804,.756}
\definecolor{verylightgrey}{rgb}{.904,.904,.856}
\usepackage[pdfauthor={V. Lebed},
pdftitle={Applications of self-distributivity to the YBE},%
colorlinks=false,linkbordercolor=lightgrey,citebordercolor=lightgrey,urlbordercolor=lightgrey]{hyperref}

\definecolor{mygreen}{rgb}{0,.455,0} 
\definecolor{myviolet}{rgb}{.45,.05,.545}
\definecolor{myred}{rgb}{.545,0,0}
\definecolor{myblue}{rgb}{.024,.15,.645}

\newcommand{\ZZ}{{\mathbb Z}}

\newcommand{\RR}{{\mathbb R}}
\newcommand{\kk}{\Bbbk}
\newcommand{\C}{{\mathcal C}}

\newcommand{\Arcs}{{\mathcal A}}
\renewcommand{\Col}{\mathpzc{Col}}

\newcommand{\Hom}{\operatorname{Hom}}
\newcommand{\Map}{\operatorname{Map}}
\newcommand{\Id}{\operatorname{Id}}
\newcommand{\Tr}{\operatorname{Tr}}
\newcommand\op{\mathrel{\triangleleft}}

\newcommand\ops{\mathrel{\triangleleft_\sigma}}

\newcommand\wdot{\mathrel{\tilde{\cdot}}}

\newcommand\Rack{{\scriptscriptstyle\mathrm{R}}}
\newcommand\Quandle{{\scriptscriptstyle\mathrm{Q}}}
\newcommand\Br{{\scriptscriptstyle\mathrm{Br}}}
\newcommand\Bir{{\scriptscriptstyle\mathrm{Bir}}}
\newcommand\Biquandle{{\scriptscriptstyle\mathrm{Biq}}}
\newcommand\RackCat{\mathbf{Rack}}
\newcommand\BirackCat{\mathbf{Birack}}
\newcommand{\Mon}{\operatorname{Mon}}
\newcommand{\Grp}{\operatorname{Grp}}

\newcommand*\Ccancel[2][black]{\renewcommand\CancelColor{\color{#1}}\cancel{#2}}

\newcounter{sarrow}
\newcommand\myleadsto[1]{
\stepcounter{sarrow}
\begin{tikzpicture}[decoration=snake]
\node (\thesarrow) {\strut#1};
\draw[->,decorate] (\thesarrow.south west) -- (\thesarrow.south east);
\end{tikzpicture}
}
\newcommand\myleads{
\begin{tikzpicture}[decoration=snake]
\draw[->,decorate] (0,0) -- (1.2,0);
\end{tikzpicture}
}

	\theoremstyle{plain}

	\theoremstyle{definition}

	\theoremstyle{remark}

\setlist{nolistsep}
\setenumerate{topsep=0pt}
\makeatother

\begin{document}

\title[Applications of self-distributivity to the YBE]{Applications of self-distributivity to Yang--Baxter operators and their cohomology}

\begin{abstract} 
Self-distributive (SD) structures form an important class of solutions to the Yang--Baxter equation, which underlie spectacular knot-theoretic applications of self-distributivity. It is less known that one can go the other way around, and construct an SD structure out of any left non-degenerate (LND) set-theoretic YBE solution. This structure captures important properties of the solution: invertibility, involutivity, biquandle-ness, the associated braid group actions. Surprisingly, the tools used to study these associated SD structures also apply to the cohomology of LND solutions, which generalizes SD cohomology. Namely, they yield an explicit isomorphism between two cohomology theories for these solutions, which until recently were studied independently. The whole story is full of open problems. One of them is the relation between the cohomologies of a YBE solution and its associated SD structure. These and related questions are covered in the present survey.
\end{abstract}



\author{Victoria Lebed}
\address{
Hamilton Mathematics Institute and
School of Mathematics,
Trinity College,
Dublin 2, Ireland}
\email{lebed.victoria@gmail.com, lebed@maths.tcd.ie}
 
\maketitle

\section{Introduction}

Self-distributive (SD) structures sporadically appeared in mathematics since the late 19th century, but became a subject of systematic study only in the 1980's, when it was realized that they give powerful and easily computable invariants of braids, knots, and their higher-dimensional analogues. Around the same time they emerged in the large cardinal theory. In these two contexts, self-distributivity can be seen as the algebraic distillation of
\begin{itemize}
\item Reidemeister move $\mathrm{III}$;
\item the properties of iterations of elementary embeddings.
\end{itemize}
A combination of these viewpoints lead to an unexpected total order on the braid groups $B_n$. Later SD reappeared in the classification of Hopf algebras and in the integration problem for Leibniz algebras. Its role in the study of the Yang--Baxter equation (YBE) was unveiled even more recently. This last application is the subject of the present survey. Our aim is to assemble relevant constructions and results scattered in different papers, and outline the state-of-the-art of the subject.

We will start with a classical introduction of SD from the knot-theoretic viewpoint (Section~\ref{S:SD_knots}). It is based on diagram colorings, of which Fox $3$-colorings are the simplest example. Besides being the most intuitive, this approach allows us to develop a braid-based graphical calculus, useful in subsequent sections. Another reason for including this standard material in our survey is to prepare the ground for analogous constructions based on the YBE. In Section~\ref{S:SD_hom} we will describe a cohomology theory for SD structures, still sticking to the knot-theoretic perspective. The content of these sections is essentially classical, so we omit most references and historical remarks for brevity; cf. the excellent surveys \cite{DehBook, CarterSurvey, QuandlesIntro, PrzSurvey, KamadaSurvey, NosakaBook} for more detail.

In Section~\ref{S:SD_YBE}, after a reminder on the Yang--Baxter equation (YBE), we will see how to construct its solutions out of SD structures; this is fairly classical. Far less known is the opposite direction: to a wide class of YBE solutions we will associate SD structures capturing their key properties. In particular, the canonical actions of the $B_n$ on the powers of a YBE solution and on the powers of its associated SD structure are isomorphic. For knot theory this is bad news: compared to SD structures, YBE solutions yield no new braid or knot invariants. However, this might change if one enriches the coloring counting invariants by Boltzmann weights, given by $2$-cocycles (or $n$-cocycles if one considers $(n-1)$-dimensional knottings in $\RR^{n+1}$). In the YBE case, these cocycles can be taken either from braided or from birack cohomology, recalled in Section~\ref{S:YBE_cohom}. These cohomologies appeared in different contexts, and both generalize SD cohomology. Reusing algebraic and graphical tools from Section~\ref{S:SD_YBE}, we will provide an explicit equivalence of the two theories in the cases where both are defined. Even though this is not a genuine application of SD to the YBE, the connecting map between the two cochain complexes would have been difficult to come up with if it has not been available from the associated SD structure study.

In spite of substantial research activity in the field, YBE solutions still retain a lot of open questions. Thus, their classification is for the moment out of reach. Also, computation techniques for the cohomology of YBE solutions are scarce, and the cohomological behaviour remains mysterious even for the most basic families of solutions. SD structures might shed new light onto these questions. For instance, it would be helpful to establish relations between the cohomology of a YBE solution and that of its associated SD structure, given that our conceptual and computational mastery of the SD cohomology is much more advanced; see the latest works \cite{FRS07,PrzYang,Nosaka15,SzymikQuillen,GarIglVen} and references therein. The existence of such relations is suggested by the explicit connecting map between the cohomologies of a YBE solution and its associated monoid; cf. Section~\ref{S:YBE_cohom} for more details. 
 The associated monoids are also at the heart of the two-step classification strategy for YBE solutions \cite{CJR_IYBG}, according to which one first classifies possible associated monoids, and then for each of them the corresponding solutions. It is natural to ask if a similar strategy with associated SD structures instead of monoids could be successful, bearing in mind the recent progress on classification of important classes of SD structures \cite{quandles,HSV}. 
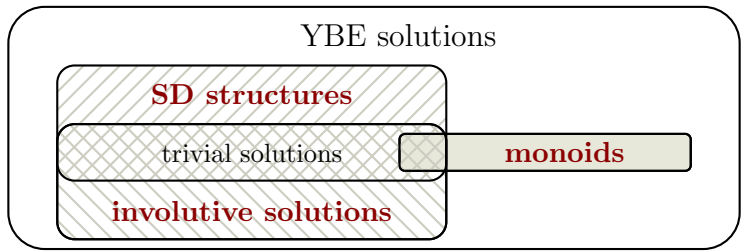
\begin{figure}[h]
\centering
\psset{unit=1.3mm}
\begin{pspicture}(70,19)(-5,-5)
\psframe[fillstyle=solid,fillcolor=verylightgrey,framearc=.3](35,3)(65,7)
\psframe[fillstyle=vlines,hatchcolor=lightgrey,framearc=.3](0,-4)(40,8)
\psframe[fillstyle=hlines,hatchcolor=lightgrey,framearc=.3](0,2)(40,14)
\psframe[framearc=.3](35,3)(65,7)
\psframe[framearc=.3](0,-4)(40,8)
\psframe[framearc=.3](-5,-5)(70,20)
\rput(20,11){\large \color{myred} \textbf{SD structures}}
\rput(20,5){trivial solutions}
\rput(20,-1){\large \color{myred} \textbf{involutive solutions}}
\rput(52,5){\large \color{myred} \textbf{monoids}}
\rput(35,17){\Large YBE solutions}
\end{pspicture}
\psset{unit=1mm}
\caption{Three major classes of YBE solutions}\label{P:Classes_YBE}
\end{figure}

The associated SD structure and monoid constructions can be seen as projection functors from the category of YBE solutions to the categories of SD structures and monoids respectively. Their sections are the corresponding inclusion functors. The ``kernel'' of the first projection is well understood: the associated SD structure is trivial precisely for involutive solutions (particularly well studied by algebraists). The ``kernel'' of the second one is still a mystery. A projection to the category of involutive solutions is also missing. To prove something about YBE solutions, it is often helpful to first look at (some of) these three particular families, summarized in Fig.~\ref{P:Classes_YBE} (here the intersection of the three families consists of the trivial $1$-element solution). A better understanding of the interactions between these ``three axes'' of the variety of YBE solutions, and in particular of the place of self-distributivity in the picture, would lead to essential progress in the area. See \cite{LebVen, LebVen3} for details and examples.

\section{Self-distributivity from a knot-theoretic viewpoint}\label{S:SD_knots}

\begin{figure}[h]
\centering
\begin{tikzpicture}[scale=0.6]
\node  at (-0.5,0)  {$a$};
\node  at (-0.5,1)  {$b$};
\node  at (1.5,0)  {$b$};
\node  at (2,1)  {$a \op b$};
\draw [thick,  rounded corners](0,0)--(0.25,0)--(0.4,0.4);
\draw [thick,  rounded corners,->](0.6,0.6)--(0.75,1)--(1,1);
\draw [thick,  rounded corners,->](0,1)--(0.25,1)--(0.75,0)--(1,0);
\end{tikzpicture} \hspace*{2cm}
\begin{tikzpicture}[scale=0.6]
\node  at (1.5,0)  {$a$};
\node  at (1.5,1)  {$b$};
\node  at (-.5,0)  {$b$};
\node  at (-1,1)  {$a \op b$};
\draw [thick,  rounded corners](0,1)--(0.25,1)--(0.4,0.6);
\draw [thick,  rounded corners,->](0.6,0.4)--(0.75,0)--(1,0);
\draw [thick,  rounded corners,->](0,0)--(0.25,0)--(0.75,1)--(1,1);
\end{tikzpicture}
\caption{Coloring rules}\label{P:Col}
\end{figure}
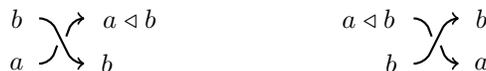

Consider a set $S$ endowed with a binary operation~$\op$. A group with the conjugation operation $a \op b = b^{-1}ab$ will be our motivating example. Taking inspiration from the Wirtinger presentation of the knot group, which is read off any of the knot's diagrams, let us consider $S$-colorings of a braid or knot diagram $D$. Concretely, we call the \emph{arcs} of $D$ its parts delimited by the crossing points, denote by $\Arcs(D)$ the set of all such arcs, and define an \emph{$S$-coloring} of $D$ to be any map $\Arcs(D) \to S$ satisfying around each crossing the rules from Fig.~\ref{P:Col}.

\begin{figure}[h]
\centering
\begin{tikzpicture}[scale=0.55]
\draw [thick,  rounded corners](0,0)--(0.25,0)--(0.4,0.4);
\draw [thick,  rounded corners](0.6,0.6)--(0.75,1)--(1.25,1)--(1.4,1.4);
\draw [thick,  rounded corners,->](1.6,1.6)--(1.75,2)--(3,2);
\draw [thick,  rounded corners](0,1)--(0.25,1)--(0.75,0)--(2.25,0)--(2.4,0.4);
\draw [thick,  rounded corners,-<](2.6,0.6)--(2.75,1)--(3,1);
\draw [thick,  rounded corners,-<](0,2)--(1.25,2)--(1.75,1)--(2.25,1)--(2.75,0)--(3,0);
\node  at (-0.5,0)  {$a$};
\node  at (-0.5,1)  {$b$};
\node  at (-0.5,2)  {$c$};
\node [right]  at (3,0)  {$c$};
\node [right]  at (3,1)  {$b \op c$};
\node [right,myblue]  at (3,2)  {$(a \op b) \op c$};
\node  at (8,1){\Large $\overset{\text{R}\mathrm{III}}{\sim}$};
\end{tikzpicture}\qquad
\begin{tikzpicture}[scale=0.55]
\draw [thick,  rounded corners](1,1)--(1.25,1)--(1.4,1.4);
\draw [thick,  rounded corners,-<](1.6,1.6)--(1.75,2)--(3.25,2)--(3.75,1)--(4,1);
\draw [thick,  rounded corners](1,0)--(2.25,0)--(2.4,0.4);
\draw [thick,  rounded corners](2.6,0.6)--(2.75,1)--(3.25,1)--(3.4,1.4);
\draw [thick,  rounded corners,-<](3.6,1.6)--(3.75,2)--(4,2);
\draw [thick,  rounded corners,->](1,2)--(1.25,2)--(1.75,1)--(2.25,1)--(2.75,0)--(4,0);
\node  at (0.5,0)  {$a$};
\node  at (0.5,1)  {$b$};
\node  at (0.5,2)  {$c$};
\node [right]  at (4,0)  {$c$};
\node [right]  at (4,1)  {$b \op c$};
\node [right,myblue]  at (4,2)  {$(a \op c) \op (b \op c)$};
\end{tikzpicture}
\caption{The R$\mathrm{III}$ move respects $(S,\op)$-colorings iff the operation $\op$ is self-distributive}\label{P:RIII}
\end{figure}
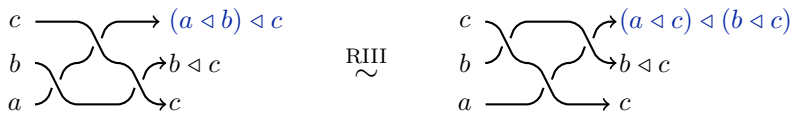

To extract invariants from such colorings, we need each relevant Reidemeister move to change colorings only locally, i.e., inside the ball where the R-move took place. As can be seen from Fig.~\ref{P:RIII}, for the R$\mathrm{III}$ move this is equivalent to the \emph{self-distributivity (SD)} of $\op$:
\begin{equation}\label{E:SD}
(a\op b)\op c = (a \op c)\op(b \op c).
\end{equation}
Note that on each side of Fig.~\ref{P:RIII}, the colors of the three leftmost arcs uniquely determine all the colors in the diagram.

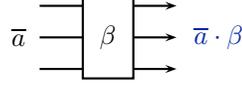
\begin{figure}[h]
\centering
\psset{unit=1.4mm}
\begin{pspicture}(17,8)(0,2)
\psline{->}(4,2)(17,2)
\psline{->}(4,5)(17,5)
\psline{->}(4,8)(17,8)
\psframe[fillstyle=solid](8,1)(13,9)
\rput(10.5,5){$\beta$}
\rput(21,5){$\color{myblue} \overline{a} \cdot \beta$}
\rput(2,5){$\overline{a}$}
\end{pspicture}
\psset{unit=1mm}
\caption{$S$-colorings define an action of (positive) $n$-strand braids on $S^n$}\label{P:action}
\end{figure}

In particular, this means that color propagation from left to right using coloring rules, as depicted in Fig.~\ref{P:action}, defines an action of the positive braid monoids $B_n^+$ on the cartesian powers $S^n$ of an SD structure $(S,\op)$. 

To extend this action to an action of the whole braid groups $B_n$, we should check the effect of the R$\mathrm{II}$ move on colorings. Finally, to work with knots (as usual, by \emph{knots} we mean knots and links in this paper), we should add the R$\mathrm{I}$ move in the picture. Algebraic counterparts of these moves are summarized in Table~\ref{T:Top_vs_Alg}. Note that the idempotence axiom can be topologically interpreted in two ways:
\begin{enumerate}
\item as the R$\mathrm{I}$ move;
\item as the triviality of the $B_n$-action on the diagonal part of $S^n$, that is, on the image of the  map $S \to S^n$, \, $a \mapsto (a,\ldots,a)$.
\end{enumerate}
 Note also that the terminology presented in the table is ``accumulative''. Thus, to be called a \emph{quandle}, the data $(S,\op)$ should satisfy the axioms from the corresponding row and all the rows above.
\begin{table}[h]
\centering
\renewcommand{\arraystretch}{1.2}
\begin{tabular}{|c|c|c|l}
\cline{1-3}
$ B_n^+ \to \operatorname{End}(S^n)$ \,&\,  R$\mathrm{III}$ \,&\, $(a\op b)\op c = (a \op c)\op(b \op c)$ & \textit{shelf} \\ \cline{1-3}
$B_n \to \operatorname{Aut}(S^n) $ \,&\, R$\mathrm{II}$ \,&\,  $\forall b,\, a \mapsto a \op b$ is invertible & \textit{rack} \\ \cline{1-3}
$S \hookrightarrow (S^n)^{B_n}$ \,&\, R$\mathrm{I}$ & $a\op a = a$ & \textit{quandle} \\
\cline{1-3}
\end{tabular}
\renewcommand{\arraystretch}{1}
\caption{A correspondence between braid actions, R-moves, and algebraic axioms}\label{T:Top_vs_Alg}
\end{table} 

Self-distributivity unifies diverse braid and knot theory gadgets in a convenient combinatorial framework. For braids, some of them are assembled in Table~\ref{T:SD_for_braids}. The corresponding $B_n^{(+)}$-actions yield the indicated representations or other braid group properties in a more or less straightforward way. Observe that for the Burau representation, modding out the diagonal part (mentioned in Table~\ref{T:Top_vs_Alg}) yields reduced Burau. For a discussion of potential braid-theoretic applications of Laver tables, see \cite{DehLeb,LaverSurvey}. The last row is the author's work in progress.
\begin{table}[h]
\centering
\begin{tabular}{|c|c|c|c|}
\hline
\rowcolor{lightgrey} $S$ & $a \op b$ & name & in braid theory \\
\hline
group & $b^{-1} ab$ & conjugation quandle & Artin: $B_n \hookrightarrow \operatorname{Aut}(F_n)$\\
\hline 
$_{\ZZ[t^{\pm 1}]}$-mod. & $ta + (1-t)b$ & Alexander quandle &  Burau: $B_n \to \operatorname{GL}_n(\ZZ[t^{\pm }])$\\
\hline
$\ZZ$ & $a+1$ & free rack & $\operatorname{lg(w)}, \operatorname{lk}_{i,j}$ \\
\hline
\multicolumn{3}{|c|}{free shelf} &  Dehornoy: order on $B_n$\\
\hline
\multicolumn{3}{|c|}{Laver table} & ???\\
\hline
\multicolumn{3}{|c|}{twisted Alexander quandle} &Lawrence--Krammer--Bigelow\\
\hline 
\end{tabular}
\caption{Examples of SD structures, with constructions in braid theory extracted from the corresponding $B_n^{(+)}$-actions}\label{T:SD_for_braids}
\end{table} 

Let us now turn to knots. Table~\ref{T:Top_vs_Alg} implies that, given a quandle $(S,\op)$ and two knot diagrams $D, D'$ related by an R-move, there is an explicit bijection
\begin{equation}\label{E:Col_Bijection}
 \Col_S(D) \overset{1:1}{\longleftrightarrow}\Col_S(D') 
\end{equation}
between their $S$-coloring sets. In particular, the number of $S$-colorings yields a knot invariant. These invariants are interesting for several reasons: they are
\begin{itemize}
\item  easy to evaluate, especially by a computer;
\item  many, as small quandles are numerous;
\item  powerful, since related to the fundamental quandle $Q(K)$ of a knot $K$, which is a weak universal knot invariant.
\end{itemize}
More precisely, the number of $S$-colorings can be computed in three different ways:  
\begin{equation}\label{E:Col_3Versions}
 \# \Col_S(D) = \# \Hom_{\mathrm{Qu}}(Q(K),S) = \Tr(\rho_S(\beta)).
\end{equation}
Here
\begin{itemize}
\item $D$ is any diagram of the knot $K$;
\item $\Hom_{\mathrm{Qu}}$ is the set of quandle homomorphisms between two quandles;
\item $\beta \in B_n$ (for some $n$) is any braid with $\operatorname{closure}(\beta) = K$;
\item $\rho_S \colon B_n \to \operatorname{Aut}(S^n)$ is the coloring action discussed above;
\item $\Tr$ is the trace (i.e., the number of fixed points) of an endomorphism of $S^n$. 
\end{itemize} 
In particular, if two quandles yield isomorphic braid group actions, then they yield the same knot invariants.

Table~\ref{T:SD_for_knots} contains basic examples of SD-interpreted knot invariants. 
\begin{table}[h]
\centering
\begin{tabular}{|c|c|}
\hline
\rowcolor{lightgrey} quandle & in knot theory \\
\hline
conjugation & representations of the knot group\\
\hline 
Alexander & Alexander polynomial\\
\hline
\end{tabular}
\caption{Examples of quandles and corresponding knot invariants}\label{T:SD_for_knots}
\end{table}

Self-distributive techniques can also be adapted to other topological objects: virtual and welded braids and knots; knotted surfaces and, more generally, $(n-1)$-dimensional knottings in $\RR^{n+1}$ \cite{PrzRos}; knotted graphs and foams \cite{CYL} etc.

\section{Self-distributive cohomology}\label{S:SD_hom}

Take a quandle $(S,\op)$ and two knot diagrams $D, D'$ related by an R-move. So far, the only invariant we have extracted from the explicit bijection~\eqref{E:Col_Bijection} between the $S$-coloring sets of the two diagrams is the cardinality of these sets. We will now show how to get much more information out of it.

Let us start with very general data: a shelf $(S,\op)$ and a map $\phi \colon S \times S \to \ZZ_n$ (or any abelian group instead of $\ZZ_n$), called the \emph{weight map}. The \emph{$\phi$-weight} of a braid or knot diagram $D$ endowed with an $S$-coloring $\C$ is defined as
\psset{unit=0.4mm,linewidth=.4,border=3pt}
\begin{equation}\label{E:Weights}
\omega_{\phi}(D,\C)=\sum_{\begin{pspicture}(10,10)(-3,-1)
\psline[linewidth=.4](0,0)(10,10)
\psline[linewidth=.4,border=0.8,arrowsize=2]{->}(0,10)(10,0)
\rput(-2,10){$\scriptstyle b$}
\rput(-2,0){$\scriptstyle a$}
\end{pspicture}} \pm\phi(a, b).
\end{equation}
Here the sum is taken over all crossings of $D$, and $\pm$ is the crossing sign. 

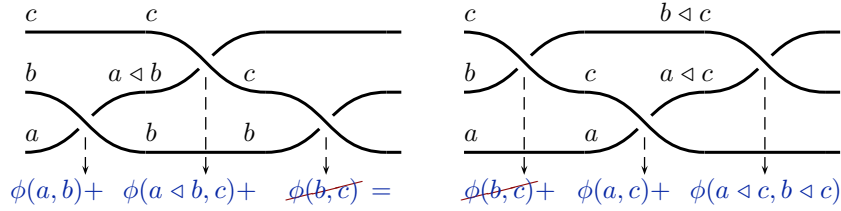
\begin{figure}[h]
\centering
\psset{unit=1mm,linewidth=.4,border=3pt}
\begin{pspicture}(60,17)(-1,-7)
\psbezier (0,0)(8,0)(8,8)(16,8)
\psbezier (0,8)(8,8)(8,0)(16,0)
\psline (0,16)(16,16)
\psbezier (16,8)(24,8)(24,16)(32,16)
\psbezier (16,16)(24,16)(24,8)(32,8)
\psline (16,0)(32,0)
\psbezier (32,0)(40,0)(40,8)(48,8)
\psbezier (32,8)(40,8)(40,0)(48,0)
\psline (32,16)(48,16)
\psline (48,0)(50,0)
\psline (48,8)(50,8)
\psline (48,16)(50,16)
\put(0,1.5){$a$}
\put(0,9.5){$b$}
\put(0,17.5){$c$}
\put(16,17.5){$c$}
\put(16,1.5){$b$}
\put(11,9.5){$a \op b$}
\put(29,9.5){$c$}
\put(29,1.5){$b$}
\psline[linewidth=0.5pt,linestyle=dashed,border=0pt]{->}(8,2)(8,-3)
\put(-2,-6){\color{myblue}$\phi(a, b)+$}
\psline[linewidth=0.5pt,linestyle=dashed,border=0pt]{->}(24,10)(24,-3)
\put(13,-6){\color{myblue}$\phi(a \op b, c)+$}
\psline[linewidth=0.5pt,linestyle=dashed,border=0pt]{->}(40,2)(40,-3)
\put(35,-6){\color{myblue}$\Ccancel[myred]{\phi(b,c)} \, =$}
\end{pspicture} \hspace*{-.4cm}
\begin{pspicture}(50,17)(0,-7)
\psbezier (0,8)(8,8)(8,16)(16,16)
\psbezier (0,16)(8,16)(8,8)(16,8)
\psline (0,0)(16,0)
\psbezier (16,0)(24,0)(24,8)(32,8)
\psbezier (16,8)(24,8)(24,0)(32,0)
\psline (16,16)(32,16)
\psbezier (32,8)(40,8)(40,16)(48,16)
\psbezier (32,16)(40,16)(40,8)(48,8)
\psline (32,0)(48,0)
\psline (48,0)(50,0)
\psline (48,8)(50,8)
\psline (48,16)(50,16)
\put(0,1.5){$a$}
\put(0,9.5){$b$}
\put(0,17.5){$c$}
\put(16,1.5){$a$}
\put(16,9.5){$c$}
\put(26,9.5){$a \op c$}
\put(26,17.5){$b \op c$}
\psline[linewidth=0.5pt,linestyle=dashed,border=0pt]{->}(8,10)(8,-3)
\put(0,-6){\color{myblue}$\Ccancel[myred]{\phi(b, c)}+$}
\psline[linewidth=0.5pt,linestyle=dashed,border=0pt]{->}(24,2)(24,-3)
\put(15,-6){\color{myblue}$\phi(a,c)+$}
\psline[linewidth=0.5pt,linestyle=dashed,border=0pt]{->}(40,10)(40,-3)
\put(30,-6){\color{myblue}$\phi(a \op c, b \op c)$} 
\end{pspicture}
\psset{border=0pt}
\caption{R$\mathrm{III}$ move and the rack $2$-cocycle condition}\label{P:Rack2Cocycle}
\end{figure}

We want the coloring bijection~\eqref{E:Col_Bijection} (valid for braids as well as knots) to induce the equality of the multi-sets of $\phi$-weights:
\begin{equation}\label{E:Col_Bijection_Weights}
\{\, \phi(D,\C)\, |\, \C \in \Col_S(D)\, \} = \{\, \phi(D',\C)\, |\, \C \in \Col_S(D')\, \}.
\end{equation}
For braids, one can restrict the colorings considered to those with prescribed values on the leftmost and rightmost ends. Instead of multi-sets, one could also work with polynomials 
\[\displaystyle \sum_{\C \in \Col_S(D)} t^{\omega_{\phi}(D,\C)} \in \ZZ[t^{\pm 1}] / (t^n-1).\] 
The desired equality holds if the relevant R-moves do not change $\omega_{\phi}(D,\C)$, which imposes algebraic conditions on $\phi$ summarized in Table~\ref{T:Top_vs_Alg_Cohom}. 
\begin{table}[h]
\centering 
\begin{tabular}{|c|c|l}
\cline{1-2}
R$\mathrm{III}$ \,&\,$\phi(a, b)+\phi(a \op b, c)=\phi(a,c)+\phi(a \op c, b \op c)$ & \textit{rack $2$-cocycle} \\ \cline{1-2}
R$\mathrm{II}$ \,&\, automatic &  \\ \cline{1-2}
R$\mathrm{I}$ & $\phi(a, a) = 0$ & \textit{quandle $2$-cocycle} \\ \cline{1-2}
\end{tabular}
\renewcommand{\arraystretch}{1}
\caption{R-moves translated as condition on the weight map}\label{T:Top_vs_Alg_Cohom}
\end{table} 
Note that a quandle $2$-cocycle should satisfy both conditions from the table. Also, the R$\mathrm{II}$ move is taken care of by the sign choice in $\pm \phi$. For the verifications for the R$\mathrm{III}$ move, see Fig.~\ref{P:Rack2Cocycle}. 

\begin{figure}[h]
\centering
\psset{unit=1mm,linewidth=.4,border=3pt}
\begin{pspicture}(0,-6)(42,11)
\psbezier (0,0)(8,0)(8,8)(16,8)
\psbezier (0,8)(8,8)(8,0)(16,0)
\psbezier (16,0)(24,0)(24,8)(32,8)
\psbezier (16,8)(24,8)(24,0)(32,0)
\put(0,1.5){$0$}
\put(0,9.5){$1$}
\put(16,1.5){$1$}
\put(16,9.5){$0$}
\put(32,1.5){$0$}
\put(32,9.5){$1$}
\psline[linewidth=0.5pt,linestyle=dashed,border=0pt]{->}(8,2)(8,-3)
\put(8,-6){\color{myblue} $1$}
\psline[linewidth=0.5pt,linestyle=dashed,border=0pt]{->}(24,2)(24,-3)
\put(24,-6){\color{myblue} $0$}
\put(42,4){$\neq$}
\end{pspicture} \hspace*{1cm}
\begin{pspicture}(0,-6)(32,11)
\psline (0,0)(32,0)
\psline (0,8)(32,8)
\put(0,1.5){$0$}
\put(0,9.5){$1$}
\put(32,1.5){$0$}
\put(32,9.5){$1$}
\put(16,-6){\color{myblue} $0$}
\end{pspicture}
\caption{Braids distinguished by a rack cocycle invariant, but not by the underlying coloring invariant}\label{P:CocycleVsCounting}
\end{figure}
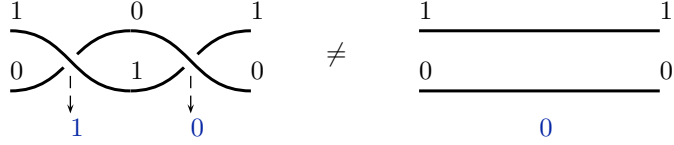

For rack/quandle $2$-cocycles $\phi$, the multi-sets of $\phi$-weights \eqref{E:Col_Bijection_Weights} are called the \emph{rack}/\emph{quandle cocycle invariants} of the braid/knot represented by $D$. They recover coloring invariants: just take the zero quandle $2$-cocycle $\phi(a,b)=0$. But for well chosen $2$-cocycles they are strictly stronger, as shown in Fig.~\ref{P:CocycleVsCounting}, where
\begin{itemize}
\item $S=\{0,1\}$ is the trivial quandle: $a \op b = a$;
\item the quandle $2$-cocycle is defined by $\phi(0,1)=1$ and $\phi(a,b)=0$ elsewhere.
\end{itemize}
Also, contrary to coloring invariants, quandle cocycle invariants detect knot chirality. See \cite{QuandleConj} for computational aspects of these invariants, which suggest that quandle cocycle invariants for finite quandles (and in fact for only a small number of those) might distinguish all knots. See \cite{IK} for an interpretation of fairly involved knot invariants, such as the complex volume and the Chern--Simons invariant, in terms of quandle cocycles.

The word ``cocycle'' above is certainly not accidental. It refers to the \emph{rack cohomology} $H^k_\Rack(S,\ZZ_n)$ of a shelf $(S,\op)$, computed using the following cochain complex:
\begin{align*}
&C^k_\Rack (S,\ZZ_n) = \Map(S^{k},\ZZ_n),\\
&(d^k_\Rack f)(a_1, \ldots, a_{k+1}) = \sum_{i=1}^{k+1} (-1)^{i-1} (f(a_1,\ldots,\widehat{a_{i}}, \ldots,a_{k+1})\\
&\hspace*{3cm} - f(a_1 {\op a_i},\ldots,a_{i-1} {\op a_i}, a_{i+1}, \ldots, a_{k+1})).
\end{align*}
Here $\widehat{a_{i}}$ means that the entry $a_{i}$ is omitted. To get the \emph{quandle cohomology} $H^k_\Quandle(S,\ZZ_n)$ of a quandle $(S,\op)$, simply restrict this complex to
\begin{align*}
&C^k_\Quandle (S,\ZZ_n) = \{\,f \colon S^{ k} \to \ZZ_n\,|\,f(\ldots,a,a,\ldots)=0  \,\}.
\end{align*}
Generalizing the above procedure, one gets invariants of $(n-1)$-dimensional knottings in $\RR^{n+1}$ out of rack/quandle $n$-cocycles \cite{PrzRos}.

Rack cohomology has several other, purely algebraic applications:
\begin{itemize}
\item pointed Hopf algebra classification \cite{AndrGr};
\item extension and deformation theories for SD structures \cite{CarterDiagrammatic,Jackson}.
\end{itemize}

\section{From self-distributivity to the Yang--Baxter equation, and back}\label{S:SD_YBE}

\begin{figure}[h]
\centering
\begin{tikzpicture}[scale=0.6]
\node  at (-0.5,0)  {$a$};
\node  at (-0.5,1)  {$b$};
\node  at (1.5,0)  {$b_a$};
\node  at (1.5,1)  {$a^b$};
\draw [thick,  rounded corners](0,0)--(0.25,0)--(0.4,0.4);
\draw [thick,  rounded corners,->](0.6,0.6)--(0.75,1)--(1,1);
\draw [thick,  rounded corners,->](0,1)--(0.25,1)--(0.75,0)--(1,0);
\end{tikzpicture} \hspace*{2cm}
\begin{tikzpicture}[scale=0.6]
\node  at (1.5,0)  {$a$};
\node  at (1.5,1)  {$b$};
\node  at (-.5,0)  {$b_a$};
\node  at (-.5,1)  {$a^b$};
\draw [thick,  rounded corners](0,1)--(0.25,1)--(0.4,0.6);
\draw [thick,  rounded corners,->](0.6,0.4)--(0.75,0)--(1,0);
\draw [thick,  rounded corners,->](0,0)--(0.25,0)--(0.75,1)--(1,1);
\end{tikzpicture}
\caption{Symmetric coloring rules}\label{P:ColYBE}
\end{figure}
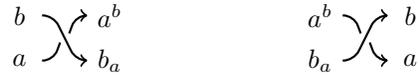

In the coloring rules from Fig.~\ref{P:Col}, we were somewhat discriminating the upper strands: contrary to the lower strands, they were not allowed to change color. One could define $S$-colorings by the more symmetric rules from Fig.~\ref{P:ColYBE} instead. 
Here we work with a set $S$ endowed with a map
\begin{align*}
\sigma \colon S^2 &\longrightarrow S^2,\\
(a,b) &\longmapsto (b_a,a^b).
\end{align*}
The SD setting is recovered by taking
\begin{equation}\label{E:Sigma_op}
\sigma_{SD}(a,b) =\sigma_{\op}(a,b) = (b,a \op b).
\end{equation}

Again, every R-move corresponds to a condition on $\sigma$, as indicated in Table~\ref{T:Top_vs_Alg_YBE}.
\begin{table}[h]
\centering 
\begin{tabular}{|c|c|c}
\cline{1-2}
R$\mathrm{III}$ \,&\, set-theoretic YBE: $\sigma_1  \sigma_2  \sigma_1 = \sigma_2  \sigma_1 \sigma_2$ & \textit{braided set} \\ \cline{1-2}
R$\mathrm{II}$ \,&\, $\sigma$ invertible \& $\forall b,\, a \mapsto a^b$ and $a \mapsto a_b$ invertible & \textit{birack} \\ \cline{1-2}
R$\mathrm{I}$ \,&\,  $\exists$ bijection $t \colon S\to S$ such that $\sigma(t(a),a) = (t(a),a)$ & \textit{biquandle} \\ \cline{1-2}
\end{tabular}
\caption{R-moves translated as conditions on the map $\sigma$}\label{T:Top_vs_Alg_YBE}
\end{table} 

As usual, the terminology on the right of the table is accumulative. Algebraists often talk about \emph{invertible non-degenerate solutions} instead of \emph{biracks}. The birack axioms can be interpreted in terms of $S$-colorings as follows: any two neighboring colors around a crossing uniquely determine the two remaining colors.

The most interesting move is R$\mathrm{III}$, translated by the set-theoretic \emph{Yang--Baxter equation (YBE)}:
\begin{equation}\label{E:YBE}
\sigma_1  \sigma_2  \sigma_1 = \sigma_2  \sigma_1 \sigma_2  \colon S^{3} \to S^{3}, \qquad \text{ where } \sigma_1 = \sigma \times  \Id_S, \, \sigma_2 = \Id_S \times \sigma.
\end{equation}
The more traditional linear YBE appears when $S$ is a vector space, all cartesian products are replaced with tensor products, and $\sigma$ is a linear map. Originating from physics, this equation is now present in purely mathematical domains as well. Quantum groups, for example, were specifically designed to produce YBE solutions.

The classification of YBE solutions is at present out of reach. It has been obtained only for $S$ of dimension~$2$, using computer-aided Gr\"{o}bner basis methods \cite{Hietarinta}. But it is still unclear how this variety of $96$ solutions (up to certain transformations) is organized. A reasonable first step, as suggested by Drinfel$'$d, could be the classification of all set-theoretic solutions, followed by an analysis of linear solutions obtained by the following procedure:

\vspace*{-.5cm}
\[\text{set-theoretic solutions } \quad\myleadsto{linearize}\quad\myleadsto{\,deform \,}\quad \text{ linear solutions}.\]

For instance, this procedure transforms the flip solution into the linear solutions originating from quantum groups:
\[\sigma(a,b) = (b,a) \quad\myleads\quad R\text{-matrices}.\]
The flip also yields a more exotic solution family \cite{Crans,CatSelfDistr,Leb}:
\[\sigma(a,b) = (b,a) \quad\myleads\quad \sigma_{Lie}(a \otimes b) = b \otimes a + \hbar 1 \otimes [a, b], \] 
where $[\,]$ is a bilinear product, $\hbar$ is an invertible scalar, and $1$ is a central element ($[1, a] = [a,1]=0$), which can be adjoined to $S$ if needed. For this map, the YBE has a particularly nice interpretation:
\[\text{YBE for } \sigma_{Lie} \qquad \Longleftrightarrow \qquad \text{ Leibniz relation for } [\,].\]
Here the Leibniz relation is taken in the form
\[ [a,[b,c]]=[[a,b],c]-[[a,c],b].\]

Another exotic example is the set-theoretic solution
\[\sigma_{Ass}(a,b) = (a \ast b,1),\] 
where $\ast$ is a binary operation on~$S$, and $1$ is a left unit ($1 \ast a = a$). Here again the YBE has a remarkable interpretation:
\[\text{YBE for } \sigma_{Ass} \qquad \Longleftrightarrow \qquad \text{ associativity for }\ast.\]

Finally, recalling the map $\sigma_{SD}$ from \eqref{E:Sigma_op}, one has
\[\text{YBE for } \sigma_{SD} \qquad \Longleftrightarrow \qquad \text{ SD for }\op.\]
This is not surprising: both conditions are dictated by the R$\mathrm{III}$ move. 

Thus, YBE solutions encapsulate such important algebraic structures as monoids, Lie (and even Leibniz) algebras, and SD structures; see Table~\ref{T:Ex_YBE}, and also Fig.~\ref{P:Classes_YBE} and the preceding discussion. For consequences of this unification, see \cite{Leb}. 
\begin{table}[h]
\centering 
\begin{tabular}{|c|c|}
\hline
\rowcolor{lightgrey} algebraic structure & YBE solution\\
\hline shelf & $\sigma_{SD}(a,b)=(b, a \op b)$\\
\hline monoid & $\sigma_{Ass}(a,b) = (a \ast b,1)$\\
\hline Lie algebra & $\sigma_{Lie}(a \otimes b) = b \otimes a + \hbar 1 \otimes [a, b]$\\
\hline
\end{tabular}
\caption{YBE solutions constructed out of basic algebraic structures}\label{T:Ex_YBE}
\end{table} 

Let us now return to low-dimensional topology. Its relations with the YBE are twofold. On the one hand, Table~\ref{T:Top_vs_Alg_YBE} implies that
\begin{itemize}
\item for a birack $(S,\op)$, $S$-colorings define a $B_n$-action on $S^n$;
\item a biquandle $(S,\op)$ yields coloring counting invariants of knots.
\end{itemize} 
In the opposite direction, braids provide a graphical calculus which is instrumental in exploring YBE solutions, as we will see in this and the next sections. In particular, we will now use it to explain why, in spite of the seemingly greater generality, birack/biquandle colorings give nothing new compared to SD structures.

To do this, we will complete the category inclusion $\RackCat \hookrightarrow \BirackCat$ given by \eqref{E:Sigma_op} with a retraction going in the opposite direction, following \cite{Sol,LYZ, LebVen}. This can be done for more general \emph{left non-degenerate} braided sets (i.e., having the maps $b \mapsto b_a$ invertible for all $a$), but for our purposes biracks are sufficient. The construction and its properties can be summarized as follows:
\begin{itemize}
\item given a birack $(S,\sigma)$, define a binary operation $\op_{\sigma}$ on $S$ by the diagram in Fig.~\ref{P:Birack_Rack}; there the colors $a$ and $b$ uniquely determine the colors of all other arcs, and $a \op_{\sigma} b$ is taken to be the rightmost bottom color;
\begin{figure}[h]
\centering
\begin{tikzpicture}[xscale=.8, yscale=.6]
\node  at (-0.3,0)  {$a$};
\node  at (1.25,-.3)  {$b$};
\node  at (3.3,0)  {$\color{myblue} a \op_{\sigma} b$};
\draw [thick,  rounded corners](0,0)--(0.25,0)--(0.4,0.4);
\draw [thick,  rounded corners, ->](0.6,0.6)--(0.75,1)--(1.75,1)--(2.25,0)--(2.5,0);
\draw [thick,  rounded corners](0,1)--(0.25,1)--(0.75,0)--(1.75,0)--(1.9,0.4);
\draw [thick,  rounded corners, ->](2.1,0.6)--(2.25,1)--(2.5,1);
\end{tikzpicture}
\caption{A rack operation on a birack}\label{P:Birack_Rack}
\end{figure}
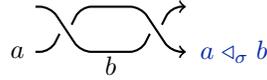
\item the resulting structure $(S,\op_{\sigma})$ is a rack, called the \emph{structure rack} or the \emph{associated rack} of $(S,\sigma)$;
\item this defines a functor $\BirackCat \twoheadrightarrow \RackCat$, which can be seen as a projection along involutive biracks, in the sense that:
\begin{itemize}
\item $\op_{\sigma_{\op}} = \op$, i.e., the composition $\RackCat \hookrightarrow \BirackCat \twoheadrightarrow \RackCat$ is the identity functor (in words, the structure rack of a birack constructed out of a rack is the original rack);
\item $\op_{\sigma}$ is trivial ($a \op_{\sigma} b = a$) $\qquad \Longleftrightarrow \qquad$ $\sigma^2= \Id$;
\end{itemize}
\item the structure rack remembers a lot about the birack:
\begin{itemize}
\item $(S,\op_{\sigma})$ is a quandle $\qquad \Longleftrightarrow \qquad$ $(S,\sigma)$ is a biquandle;
\item $\sigma$ and $\op_{\sigma}$ induce isomorphic $B_n$-actions on $S^n$.
\end{itemize}
\end{itemize}
These properties imply that rack and birack colorings yield the same braid invariants. Recalling~\eqref{E:Col_3Versions}, one concludes that quandle and biquandle colorings yield the same knot invariants.

Probably the most enlightening proofs of the properties of $\op_{\sigma}$ use diagrammatic arguments. Here we will give two of them. First, the self-distributivity of $\op_{\sigma}$ can be established by following Fig.~\ref{P:Birack_Rack_Proof}. In the left diagram, the colors $a$, $b$ and $c$ (in black) uniquely determine the colors of all other arcs. 
\begin{figure}[h]
\centering
\labellist
\small\hair 2pt
\hspace*{5cm}

\pinlabel ${\color{myblue} \left( a\ops b \right)\ops c }$ at 205 565
\pinlabel $c$ at 200 450
\pinlabel ${\color{myblue} a \ops b}$ at 185 360
\pinlabel $b$ at 200 263
\pinlabel $a$ at 162 149

\pinlabel ${\color{myblue} \left( a \ops b \right)\ops c }$ at 410 565
\pinlabel ${\color{myblue} b \ops c}$ at 415 450
\pinlabel $c$ at 395 360
\pinlabel $b$ at 389 263
\pinlabel $a$ at 359 149

\pinlabel ${\color{blue} \left( a\ops b \right)\ops c =}$ at 610 565
\pinlabel ${\color{myred} \left( a\ops c \right)\ops \left( b \ops c \right)}$ at 640 535
\pinlabel ${\color{myblue} b \ops c}$ at 610 450
\pinlabel ${\color{myblue} a \ops c}$ at 575 360
\pinlabel $c$ at 590 263
\pinlabel $a$ at 555 149

\endlabellist
\centering 
\vspace*{-.3cm}
\includegraphics[scale=0.4]{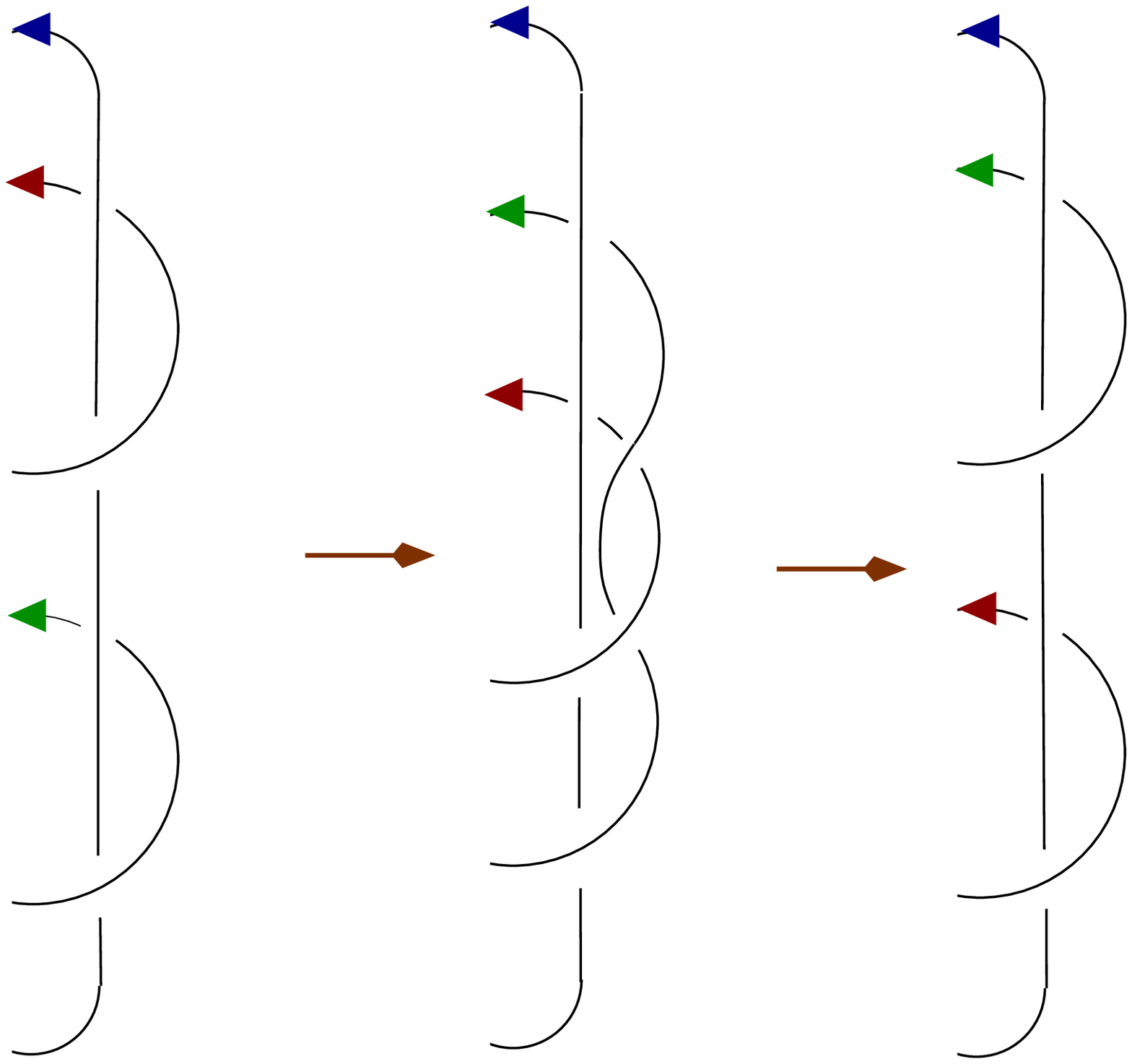}
\vspace*{-2.5cm}
\caption{A self-distributivity proof for the induced rack operation $\op_{\sigma}$}\label{P:Birack_Rack_Proof}
\end{figure}
We then move the lower semicircle through the upper one, and determine the colors of the arcs created on the way. In the end, the upper color can be determined in two ways---from the final diagram (where one recognizes the situation from the definition of $\op_{\sigma}$), or carried from the preceding one. The equality of the two expressions is the desired SD property.

\begin{figure}[h]
\centering
\begin{tikzpicture}[xscale=0.8,>=latex]
 \draw [line width=4pt,white, rounded corners=10] (4,0) -- (5,0.5) -- (0,3);   
 \draw [->, rounded corners=10] (4,0) -- (5,0.5) -- (0,3);  
 \draw [line width=4pt,white, rounded corners=10] (3,0) -- (5,1) -- (1,3); 
 \draw [->, rounded corners=10] (3,0) -- (5,1) -- (1,3); 
 \draw [line width=4pt,white, rounded corners=10] (2,0) -- (5,1.5) -- (2,3); 
 \draw [->, rounded corners=10] (2,0) -- (5,1.5) -- (2,3);
 \draw [line width=4pt,white, rounded corners=10] (1,0) -- (5,2) -- (3,3); 
 \draw [->, rounded corners=10] (1,0) -- (5,2) -- (3,3); 
 \draw [line width=4pt,white, rounded corners=10] (0,0) -- (5,2.5) -- (4,3); 
 \draw [->, rounded corners=10] (0,0) -- (5,2.5) -- (4,3);
 \node at (0,-0.4) [above] {$\scriptstyle{a_n}$};  
 \node at (2,-0.4) [above] {$\scriptstyle{\cdots}$}; 
 \node at (4,-0.4) [above] {$\scriptstyle{a_1}$}; 
 \node at (4.7,2.5) [right] {$\scriptstyle{J_n(\overline{a})}$};   
 \node at (4.7,1) [right] {$\scriptstyle{J_2(\overline{a})}$}; 
 \node at (5,1.75) [right] {$\;\scriptstyle{\vdots}$};  
 \node at (4.7,0.5) [right] {$\scriptstyle{J_1(\overline{a})}$};  
\end{tikzpicture}  \quad
\begin{tikzpicture}[xscale=0.8,>=latex]
 \draw [line width=4pt,white, rounded corners=10] (4,0) -- (5,0.5) -- (0,3);   
 \draw [->, rounded corners=10] (4,0) -- (5,0.5) -- (0,3);  
 \draw [line width=4pt,white, rounded corners=10] (3,0) -- (5,1) -- (1,3); 
 \draw [->, rounded corners=10] (3,0) -- (5,1) -- (1,3); 
 \draw [line width=4pt,white, rounded corners=10] (2,0) -- (5,1.5) -- (2,3); 
 \draw [->, rounded corners=10] (2,0) -- (5,1.5) -- (2,3);
 \draw [line width=4pt,white, rounded corners=10] (1,0) -- (5,0.6)  -- (5.5,1.1) -- (5,2) -- (3,3); 
 \draw [->, rounded corners=10,myred] (1,0) -- (5,0.6)  -- (5.5,1.1) -- (5,2) -- (3,3); 
 \draw [line width=4pt,white, rounded corners=10] (0,0) -- (5,2.5) -- (4,3);
 \draw [->, rounded corners=10] (0,0) -- (5,2.5) -- (4,3);
 \node at (1,-0.4) [above] {$\scriptstyle{a_4}$};
 \node at (2,-0.4) [above] {$\scriptstyle{a_3}$};  
 \node at (3,-0.4) [above] {$\scriptstyle{a_2}$}; 
 \node at (4,-0.4) [above] {$\scriptstyle{a_1}$}; 
 \node at (4.9,2) [right] {$\scriptstyle{J_4(\overline{a})=}$}; 
 \node at (5.3,1.6) [right] {$\scriptstyle{(a_4)_{a_3 a_2 a_1}}$};    
\end{tikzpicture} 
\caption{The guitar map}\label{P:Guitar}
\end{figure}
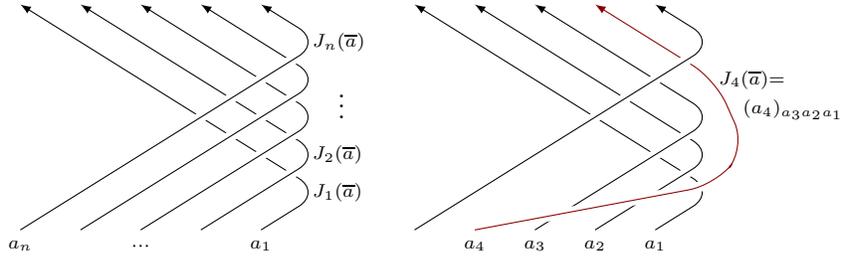

Second, to compare the $B_n$-actions induced by $\sigma$ \& $\op_{\sigma}$, consider the \emph{guitar map}
\begin{align*}
J\colon S^{n} &\longrightarrow S^{n}, \\
(a_1, \ldots, a_n) &\longmapsto (a_1, (a_{2})_{a_1}, (a_3)_{a_{2} a_1}, \ldots ). 
\end{align*}
We used simplified notations $(a_3)_{a_{2} a_1} = ((a_3)_{a_2})_{a_1}$ etc. Fig.~\ref{P:Guitar} presents $J$ graphically: when the colors $a_n, \ldots, a_1$ are propagated along the left diagram, the $n$ rightmost colors yield $J(a_1, \ldots, a_n)$, as shown in the right diagram. The guitar map appeared under different names and at different levels of generality in \cite{ESS,Sol,LYZ,LebVen}. We will meet it again in the next section. Note that it is well defined for any braided set.

To get a feeling of how this map works, let us look at it in particular cases:
\begin{enumerate}
\item $\sigma_{Ass}(a,b)=(a\ast b,1) \quad \leadsto \quad J(a,b,c)=(a, a\ast b, a\ast b\ast c)$;

\noindent this map relates two forms of the bar complex for the monoid $(S,\ast,1)$;

\item $\widetilde{\sigma_{SD}}(a,b)=(b \op a, a) \quad \leadsto \quad J(a,b,c)=(a, b\op a, (c \op b) \op a)$;

\noindent this is the \textit{remarkable map} from \cite{Prz1}; note that we are using here the mirror version of the YB operator $\sigma_{SD}$ from~\eqref{E:Sigma_op};

\item $\sigma^2=\Id \quad \leadsto \quad$ the map $\Omega$ from the right-cyclic calculus \cite{DehCycleSet}.
\end{enumerate}

Among the remarkable properties of the guitar map for a birack $(S,\sigma)$, the two that are relevant to us are:
\begin{itemize}
\item $J$ is invertible;
\item $J\sigma_i = \sigma'_i J$, where 
\[\sigma'=\sigma_{\op_{\sigma}} \colon (a,b) \mapsto (b, a \op_{\sigma} b)\]
is the YB operator built from the structure rack of $(S,\sigma)$, and
\[\sigma_i = \Id_S^{i-1} \times \sigma \times \Id_S^{n-i-1},\]
and similarly for $\sigma'_i$.
\end{itemize}
These properties imply that $\sigma$ and $\op_{\sigma}$ yield isomorphic $B_n$-actions on $S^n$, without the biracks $(S, \sigma)$ and $(S, \sigma')$ being isomorphic in general.

\begin{figure}[h]
\centering
\begin{tikzpicture}[xscale=0.7,yscale=0.9,>=latex]
 \draw [line width=4pt,white, rounded corners=10] (4,0) -- (5,0.5) -- (0,3);   
 \draw [->, rounded corners=10] (4,-0.5) -- (4,0) -- (5,0.5) -- (0,3);  
 \draw [line width=4pt,white, rounded corners=10] (2,0) -- (5,1.5) -- (2,3); 
 \draw [->, rounded corners=10] (1,-0.5) -- (2,0) -- (5,1.5) -- (2,3); 
 \draw [line width=4pt,white, rounded corners=10] (2,-0.5) -- (1,0) -- (5,2) -- (3,3); 
 \draw [->, rounded corners=10] (2,-0.5) -- (1,0) -- (5,2) -- (3,3);
 \draw [line width=4pt,white, rounded corners=10]  (0,-0.5) -- (0,0) -- (5,2.5) -- (4,3);
 \draw [->, rounded corners=10] (0,-0.5) -- (0,0) -- (5,2.5) -- (4,3);
 \draw [line width=4pt,white, rounded corners=10] (1,-0.5) -- (2,0); 
 \draw [rounded corners=10] (1,-0.5) -- (2,0); 
 \draw [myviolet, dashed] (-1,0)--(5,0); 
 \node at (-1,-0.5) [myviolet] {${\sigma_2}$};  
 \node at (-1,1.5) [myviolet]  {${J}$}; 
 \node at (0,-0.9) [above] {$\scriptstyle{a_4}$};  
 \node at (1,-0.9) [above] {$\scriptstyle{a_3}$}; 
 \node at (2,-0.9) [above] {$\scriptstyle{a_2}$}; 
 \node at (4,-0.9) [above] {$\scriptstyle{a_1}$};
 \node at (4.7,2.5) [right] {$\scriptstyle{\color{myblue}J_4(\sigma_2(\overline{a}))}$};   
 \node at (4.7,2) [right] {$\scriptstyle{\color{myblue}J_3(\sigma_2(\overline{a}))}$};   
 \node at (4.7,1.5) [right] {$\scriptstyle{\color{myblue}J_2(\sigma_2(\overline{a}))}$};  
 \node at (4.7,0.5) [right] {$\scriptstyle{\color{myblue}J_1(\sigma_2(\overline{a}))}$};  
\end{tikzpicture} 
\hspace*{20pt}
\begin{tikzpicture}[xscale=0.7,yscale=0.9,>=latex]
 \draw [->, rounded corners=10] (3,2.5) -- (2,3);  
 \draw [line width=4pt,white, rounded corners=10] (4,0) -- (5,0.5) -- (0,3);   
 \draw [->, rounded corners=10] (4,-0.5) -- (4,0) -- (5,0.5) -- (0,3);  
 \draw [line width=4pt,white, rounded corners=10] (2,-0.5) -- (5,1) -- (2,2.5) -- (3,3);
 \draw [->, rounded corners=10] (2,-0.5) -- (5,1) -- (2,2.5) -- (3,3);
 \draw [line width=4pt,white, rounded corners=10] (1,-0.5) --(2,0) -- (5,1.5) -- (3,2.5); 
 \draw [rounded corners=10] (1,-0.5) -- (2,0) -- (5,1.5) -- (3,2.5);  
 \draw [line width=4pt,white, rounded corners=10] (0,-0.5) -- (0,0) -- (5,2.5) -- (4,3);
 \draw [->, rounded corners=10] (0,-0.5) -- (0,0) -- (5,2.5) -- (4,3);
 \draw [myviolet, dashed] (-0.5,2.5)--(5.5,2.5); 
 \node at (-0.5,0.5) [myviolet]  {${J}$}; 
 \node at (0,-0.9) [above] {$\scriptstyle{a_4}$};  
 \node at (1,-0.9) [above] {$\scriptstyle{a_3}$}; 
 \node at (2,-0.9) [above] {$\scriptstyle{a_2}$}; 
 \node at (4,-0.9) [above] {$\scriptstyle{a_1}$};
 \node at (4.7,2.25) [right] {$\scriptstyle{\color{red} J_4(\overline{a})}$};   
 \node at (4.7,1.5) [right] {$\scriptstyle{\color{red} J_3(\overline{a})}$};  
 \node at (4.7,1) [right] {$\scriptstyle{\color{red} J_2(\overline{a})}$};   
 \node at (4.7,0.5) [right] {$\scriptstyle{\color{red} J_1(\overline{a})}$};   
\end{tikzpicture}
\begin{tikzpicture}[xscale=0.7,yscale=0.9,>=latex]
 \draw [latex->, rounded corners=10,mygreen] (-2,3.5) -- (-1.5,3)node[below left]
{R$\mathrm{III}$} -- (-0.5,2.5);
 \draw [latex->, rounded corners=10,mygreen] (11.5,3.5) -- (11,3)node[below right]
{R$\mathrm{III}$} -- (10,2.5);
 \draw [line width=4pt,white, rounded corners=10] (4,-0.5) -- (4,0) -- (5,0.5) -- (0,3);   
 \draw [->, rounded corners=10] (4,-0.5) -- (4,0) -- (5,0.5) -- (0,3); 
 \draw [line width=4pt,white, rounded corners=10] (2,-0.5) -- (5,1) -- (4,1.5) -- (5,2) -- (3,3); 
 \draw [->, rounded corners=10] (2,-0.5) -- (5,1) -- (4,1.5) -- (5,2) -- (3,3);
 \draw [line width=4pt,white, rounded corners=10] (1,-0.5) -- (2,0) -- (5,1.5) -- (2,3); 
 \draw [->, rounded corners=10] (1,-0.5) -- (2,0) -- (5,1.5) -- (2,3);  
 \draw [line width=4pt,white, rounded corners=10] (0,-0.5) -- (0,0) -- (5,2.5) -- (4,3);
 \draw [->, rounded corners=10] (0,-0.5) -- (0,0) -- (5,2.5) -- (4,3); 
 \draw [line width=4pt,white, rounded corners=10] (4.4,1.7) -- (4.6,1.8);
 \draw [rounded corners=10] (4.4,1.7) -- (4.6,1.8);
 \node at (0,-0.9) [above] {$\scriptstyle{a_4}$};  
 \node at (1,-0.9) [above] {$\scriptstyle{a_3}$}; 
 \node at (2,-0.9) [above] {$\scriptstyle{a_2}$}; 
 \node at (4,-0.9) [above] {$\scriptstyle{a_1}$};
 \node at (4.7,2.5) [right] {$\scriptstyle{{\color{myblue}J_4(\sigma_2(\overline{a}))} = \color{red} J_4(\overline{a})}$};   
 \node at (4.7,2) [right] {$\scriptstyle{{\color{myblue}J_3(\sigma_2(\overline{a}))} = \color{red} J_2(\overline{a}) \ops  J_3(\overline{a})}$};   
 \node at (4.7,1.5) [right] {$\scriptstyle{{\color{myblue}J_2(\sigma_2(\overline{a}))} = \color{red} J_3(\overline{a})}$};  
 \node at (4.7,1) [right] {$\scriptstyle{\color{red} J_2(\overline{a})}$};   
 \node at (4.7,0.5) [right] {$\scriptstyle{{\color{myblue}J_1(\sigma_2(\overline{a}))} = \color{red} J_1(\overline{a})}$};   
\end{tikzpicture}
\caption{A proof of the entwining relation $J \sigma_i = \sigma'_i J$ (here $n=4$, $i=2$)}\label{P:Same_Bn_actions}
\end{figure}
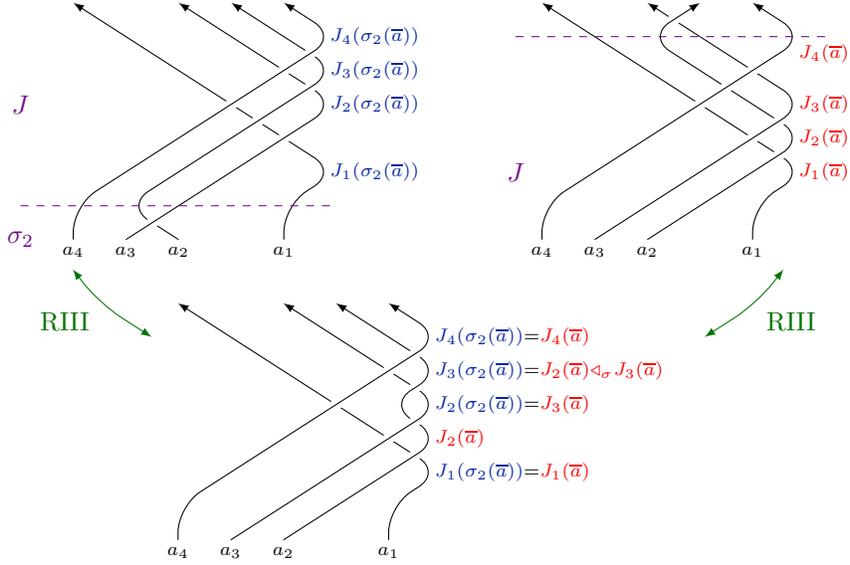

As usual, the relation $J\sigma_i = \sigma'_i J$ can be proved graphically. In Fig.~\ref{P:Same_Bn_actions}, the colors of the rightmost arcs in the bottom diagram are calculated in two ways: from the upper left (blue labels) and the upper right diagrams (red labels). In both cases, the lower colors $a_1,\ldots,a_n$ are propagated throughout the diagrams.

\section{Braided cohomology vs. birack cohomology}\label{S:YBE_cohom}

We have seen that, as far as coloring invariants are concerned, general biracks and biquandles are no better than racks and quandles. Things might change if weights are brought into the picture. Even if these upgraded invariants are not well studied yet, the cohomology theory of YBE solutions, which controls the weights, is of interest in its own right. In this section we will review its two forms, and establish their isomorphism using the guitar map from the previous section.

\begin{figure}[h]
\centering
\begin{tikzpicture}[xscale=0.5,yscale=0.3,rotate=-90] 
 \draw [myred,rounded corners] (-2,-2) -- (-2,-3) -- (4,-6) -- (4,-7);
 \draw [white,line width=4pt] (0,-2.5) -- (0,-7);
 \draw [white,line width=4pt] (3,-2.5) -- (3,-7);
 \draw (0,-2.5) -- (0,-7);
 \draw (3,-2.5) -- (3,-7);
 \draw (5,-2.5) -- (5,-7);
 \draw (8,-2.5) -- (8,-7);
 \draw (-1,-2.5) rectangle (9,-1);
 \node at (4,-1.75) {$f$};
 \node at (-2,-2) {\color{myred} $\bullet$};
 \node at (4,-11) {$ =\sum(-1)^{i-1} \Big($};
 \node at (0,-7.7) {$\scriptstyle a_{n+1}$};
 \node at (.6,-3.3) {$\color{myred} \scriptstyle a'_{n+1}$}; 
 \node at (3.6,-3.3) {$\color{myred} \scriptstyle a'_{i+1}$}; 
 \node at (1.5,-7.7) {$\scriptstyle \ldots$};
 \node at (3,-7.7) {$\scriptstyle a_{i+1}$};
 \node at (4,-7.7) {$\scriptstyle a_i$};
 \node at (6.5,-7.7) {$\scriptstyle\ldots$};
 \node at (8.5,-7.7) {$\scriptstyle a_1$};
\end{tikzpicture}\hspace*{.5cm}
\begin{tikzpicture}[xscale=0.5,yscale=0.3,rotate=-90] 
 \draw (0,-2.5) -- (0,-7);
 \draw (3,-2.5) -- (3,-7);
 \draw (5,-2.5) -- (5,-7);
 \draw (8,-2.5) -- (8,-7);
 \draw [white,line width=4pt,rounded corners] (10,-2) -- (10,-3) -- (4,-6) -- (4,-7);
 \draw [myred,rounded corners] (10,-2) -- (10,-3) -- (4,-6) -- (4,-7);
 \draw (-1,-2.5) rectangle (9,-1);
 \node at (4,-1.75) {$f$};
 \node at (10,-2) {\color{myred} $\bullet$};
 \node at (5.5,-9.5) {$-$};
 \node at (5.5,-.2) {$ \Big)$};
 \node at (0,-7.7) {$\scriptstyle a_{n+1}$};
 \node at (8.6,-3) {$\color{myred} \scriptstyle a'_{1}$}; 
 \node at (5.6,-3.3) {$\color{myred} \scriptstyle a'_{i-1}$}; 
 \node at (1.5,-7.7) {$\scriptstyle \ldots$};
 \node at (5,-7.7) {$\scriptstyle a_{i-1}$};
 \node at (4,-7.7) {$\scriptstyle a_i$};
 \node at (6.5,-7.7) {$\scriptstyle\ldots$};
 \node at (8.5,-7.7) {$\scriptstyle a_1$};
\end{tikzpicture}
\caption{A diagrammatic formula for $(d^k_\Br f)(a_1, \ldots, a_{k+1})$}\label{P:BrCohom}
\end{figure}
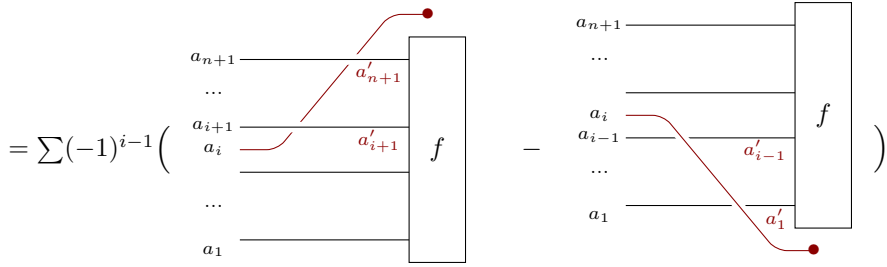

The \emph{braided cohomology} $H^k_\Br(S,\ZZ_n)$ of a braided set $(S,\sigma)$ is defined as the cohomology of the following complex \cite{CES_BrHom,Leb}:
\begin{align*}
&C^k_\Br (S,\ZZ_n) = \Map(S^{k},\ZZ_n),\\
&(d^k_\Br f)(a_1, \ldots, a_{k+1}) = \sum_{i=1}^{k+1} (-1)^{i-1} (f(a_1,\ldots,a_{i-1}, {({a_{i+1}}, \ldots, a_{k+1})_{a_i}})\\
&\hspace*{3cm} - f({(a_1,\ldots, a_{i-1})^{a_i}}, a_{i+1},\ldots,a_{k+1})).
\end{align*}
Here we used the inductive definition
\begin{align*}
(a, b, \ldots, v)_{w} &= (a_{w}, (b, \ldots, v)_{w^a}), \\ 
(a, \ldots, u, v)^{w} &= ((a, \ldots, u)^{w_v},v^{w}).
\end{align*}
The formula for the differential is best understood graphically, see Fig.~\ref{P:BrCohom}. 
Here as usual colors are propagated from left to right, and an $f$-labeled box stands for the evaluation of $f$ on the entries given by the colors of the incoming arcs, read from bottom to top. The result is an alternating sum of these evaluations. The verification of the differential property $d^{k+1}_\Br d^k_\Br=0$ now boils down to simple diagram manipulations.

For different purposes, alternative definitions of the braided complex might be more adequate. They are based on:
\begin{enumerate}
\item a cubical classifying space \cite{CES_BrHom};
\item the quantum shuffle coproduct \cite{Leb};
\item a differential graded bialgebra encoding the complex \cite{FarinatiGalofre}.
\end{enumerate}

To get the \emph{braided biquandle cohomology} of a biquandle $(S,\sigma)$, one restricts this complex to a certain subcomplex \cite{LebVen}. In degree~$2$, it is defined as
\begin{align*}
&\{\,f \colon S^{2} \to \ZZ_n\,|\,f(t(a),a)=0  \,\}.
\end{align*}

Imitating rack and quandle cocycle invariants, one extracts invariants of
\begin{itemize}
\item braids out of a birack and its braided $2$-cocycle;
\item knots out of a biquandle and its braided biquandle $2$-cocycle.
\end{itemize}
They refine the coloring invariants based on the same birack/biquandle. This refinement is trivial if it uses a $2$-coboundary. As usual, braided (biquandle) $n$-cocycles yield invariants of $(n-1)$-dimensional knottings in $\RR^{n+1}$. Thus braided cohomology offers powerful tools for braid and knot classification.

Again, one can go the other way around and apply the braid-based graphical calculus developed above to unveiling the structure of braided cohomology. For instance, this calculus renders the cup product completely intuitive:
\[\smile \, \colon C^k_\Br \otimes C^n_\Br \to C^{k+n}_\Br,\]
\begin{center}
\begin{tikzpicture}[xscale=0.6,yscale=0.3,rotate=-90] 
 \draw [myblue] (0,0) -- (-1,3);
 \draw [myblue] (2,0) -- (0,3); 
 \draw [myblue] (3,0) -- (1,3); 
 \draw [myblue] (6,0) -- (2,3); 
 \draw [white,line width=4pt] (-1,0) -- (4,3);
 \draw [white,line width=4pt] (1,0) -- (5,3);
 \draw [white,line width=4pt] (4,0) -- (6,3); 
 \draw [white,line width=4pt] (5,0) -- (7,3); 
 \draw [white,line width=4pt] (7,0) -- (8,3);
 \draw [myred] (-1,0) -- (4,3);
 \draw [myred] (1,0) -- (5,3);
 \draw [myred] (4,0) -- (6,3); 
 \draw [myred] (5,0) -- (7,3); 
 \draw [myred] (7,0) -- (8,3);
 \draw (-1.5,3) rectangle (2.5,4);
 \node at (.5,3.5) {$g$};
 \draw (3.5,3) rectangle (8.5,4);
 \node at (6,3.5) {$f$};
 \node at (-1,-.7) {$a_{k+n}$};  
 \node at (7,-.5) {$a_{1}$};  
 \node at (3,-.5) {$\cdot$};  
 \node at (2,-.5) {$\cdot$};  
 \node at (4,-.5) {$\cdot$};  
 \node at (3,-6) {$f \smile g(a_1,\ldots, a_{k+n})= \displaystyle \sum_{\text{splittings}} (-1)^{\raisebox{.1cm}{\#}  \begin{tikzpicture}[scale=0.3] 
 \draw [myblue] (0,0) -- (1,1);
 \draw [white,line width=2pt] (0,1) -- (1,0);
 \draw [myred] (0,1) -- (1,0);
\end{tikzpicture} }$};  
\end{tikzpicture} 
\end{center}
Simple diagram manipulations yield the properties of this product:
\begin{itemize}
\item $(C^*_\Br, \smile)$ is a differential graded associative algebra, graded commutative up to an explicit homotopy (which can be defined graphically);
\item $(H^*_\Br, \smile)$ is a graded commutative associative algebra; 
 here we abusively denote by $\smile$ the induced product on cohomology.
\end{itemize}
This cup product was discovered in different forms and at different levels of generality in  \cite{Clauwens,Covez,FarinatiGalofre,LebedIdempot}. When $\sigma=\sigma_{\op}$ comes from an SD operation on $S$, the commutative cup product in cohomology can be refined into a Zinbiel structure \cite{Covez2}.

Also, the diagrammatic interpretation of braided cohomology makes obvious
\begin{itemize}
\item its generalization to YBE solutions in any preadditive monoidal category (which includes set-theoretic and linear solutions);
\item its dual, homological version;
\item its functoriality;
\item its enhancement by introducing coefficients (used to color the regions of a diagram together with the arcs; this is called \emph{shadow coloring} in SD theory).
\end{itemize}

Braid and knot invariants are far from being the only application of braided cohomology. We will now outline several others.

First, braided cohomology in degree $2$ controls \emph{diagonal deformations} of a braided set $(S,\sigma)$. Namely, the map
\[\sigma_q(a,b) = q^{\phi(a,b)} \sigma(a,b),\]
where $q$ lies in the base field $\kk$ and $q^n=1$, is a YB operator on $\kk S$ if and only if $\phi$ is a $2$-cocycle \cite{FY}. Moreover, if $\phi$ is a $2$-coboundary, then this YBE solution is isomorphic to (the linearization of) the original one. It is much more challenging to describe general deformations of braided sets. A cohomology theory for that was proposed in \cite{Eisermann}. It contains braided cohomology. Because of its generality, its computation is out of reach at the moment. A better understanding of that theory would be a breakthrough in the classification of YBE solutions.

Second, we know how to construct a YBE solution (set-theoretic or linear) out of a monoid, Lie algebra, or shelf; cf. Table~\ref{T:Ex_YBE}. The braided cohomology of these solutions turns out to contain the classical cohomologies of the original structures \cite{Leb}. This unification allows the transport of constructions and results from well-studied settings to less explored ones. For instance, historically the development of SD cohomology was motivated by topological and Hopf-algebraic applications, but relied heavily on tools borrowed from the already classical cohomology theory for associative structures. However, as pointed out in \cite{Prz1}, there was no conceptual explanation of the success of this borrowing. The unifying braided cohomology setting offered such an explanation. It also guided the development of cohomology theories for new structures, such as cycle-sets and braces \cite{LebVen2}. The latter were designed in \cite{Rump_CSets,Rump_Braces} to study involutive YB operators.

Finally, for a braided set $(S,\sigma)$ with $\sigma$ involutive or idempotent, braided cohomology computes the cohomology of its \emph{structure}, or \emph{associated, monoid} 
\[\Mon({S,\sigma}) = \langle \; S \; | \; a \cdot b = b_a \cdot a^b \text{ for all } a,b \in S \; \rangle\]
\cite{FarinatiGalofre,LebedIdempot}. We used our usual notation $\sigma(a,b) = (b_a \, ,\, a^b)$. That is, if one manages to represent a monoid of interest as the structure monoid of a nice braided set, then one gets information about the cohomology of this monoid. This was applied to the cohomology of factorized monoids \cite{LebedIdempot} and plactic monoids \cite{LebedPlactic}. Note that for general $\sigma$, one still has a map from the cohomology of $\Mon({S,\sigma})$ to that of $(S,\sigma)$. It is called the \emph{quantum symmetrizer}, and can be defined graphically. Describing its kernel and image, and more generally understanding the relations between the two cohomologies, is an open question, raised in \cite{FarinatiGalofre,YangGraphYBE}.

Structure monoids are interesting for several reasons. On the one hand, they introduce group-theoretic methods into the study of the YBE, just as structure racks bring SD methods. On the other hand, especially for involutive $\sigma$, these monoids and the corresponding groups $\Grp({S,\sigma})$ and algebras $\kk \Mon({S,\sigma})$ boast interesting algebraic and geometric properties, and have already served as remarkable examples and counter-examples in different situations. In our favorite situations,
\begin{enumerate}
\item for a monoid $S$, $\Mon({S,\sigma})\cong S$;
\item for a rack $S$, $\Grp({S,\sigma_{SD}}) = \langle \; S \; | \; a\op b = b^{-1} a b\;\rangle$ is the {associated group} of the rack, which is a classical construction in the SD theory.
\end{enumerate}
Note that $\Grp({S,\sigma})$-modules serve as coefficients for an enhanced version of braided cohomology. Also, structure racks turn out to be useful in the study of structure groups and monoids \cite{LebVen3}. 

We will next describe another cohomology theory for biracks. Its advantages are
\begin{itemize}
\item more manageable formulas;
\item a natural subcomplex capturing the additional quandle axiom.
\end{itemize}
On the flip side, this theory does not generalize to other types of YBE solutions.

\begin{figure}[h]
\centering
\begin{tikzpicture}[scale=0.6]
\node  at (-0.5,0)  {$a$};
\node  at (-1,1)  {$\color{myblue} a \cdot b$};
\node  at (1.5,0)  {$b$};
\node  at (2,1)  {$\color{myblue} a \wdot b$};
\draw [thick,  rounded corners](0,0)--(0.25,0)--(0.4,0.4);
\draw [thick,  rounded corners, ->](0.6,0.6)--(0.75,1)--(1,1);
\draw [thick,  rounded corners, ->](0,1)--(0.25,1)--(0.75,0)--(1,0);
\end{tikzpicture}
\caption{Sideways operations for a birack}\label{P:Sideways}
\end{figure}
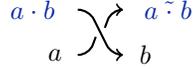

Recall that for birack colorings, any two neighboring colors around a crossing uniquely determine the two remaining colors. This justifies the existence of the \emph{sideways operations} $\cdot$ and $\wdot$, defined in Fig.~\ref{P:Sideways}.

Using these operations, the \emph{birack cohomology} $H^k_\Bir(S,\ZZ_n)$ of a birack $(S,\sigma)$ is defined as the cohomology of the following complex \cite{FRS_Species,BirackHom}:
\begin{align*}
&C^k_\Bir (S,\ZZ_n) = \Map(S^{k},\ZZ_n),\\
&(d^k_\Bir f)(a_1, \ldots, a_{k+1}) = \sum_{i=1}^{k+1} (-1)^{i-1} (f(a_1,\ldots,{\widehat{a_{i}}}, \ldots,a_{k+1})\\
&\hspace*{3cm} - f({a_i \wdot  a_1},\ldots,{a_i \wdot  a_{i-1}}, {a_i \cdot a_{i+1}},\ldots,{a_i \cdot a_{k+1}})).
\end{align*}

The \emph{biquandle cohomology} $H^k_\Biquandle(S,\ZZ_n)$ of a biquandle is the cohomology of the subcomplex
\begin{align*}
&C^k_\Biquandle (S,\ZZ_n) = \{\,f \colon S^{ k} \to \ZZ_n\,|\,f(\ldots,a,a,\ldots)=0  \,\}.
\end{align*}
In fact, there is an explicit splitting of the cochain complex \cite{LebVen} 
\[C^{\bullet}_\Bir \cong C^{\bullet}_\Biquandle \oplus C^{\bullet}_D,\]
that generalizes the splitting for the rack cohomology of quandles \cite{LN}. It would be interesting to know if the biquandle cohomology part $H^{\bullet}_\Biquandle$ completely determines the degenerate part $H^{\bullet}_D$, as is the case for quandles \cite{PrzPu}.

\psset{unit=0.4mm,linewidth=.4,border=3pt}
Now, changing the definition of $\phi$-weights from \eqref{E:Weights} to
\[\omega^*_{\phi}(D,\C)=\sum_{\begin{pspicture}(10,10)(-3,-1)
\psline[linewidth=.4](0,0)(10,10)
\psline[linewidth=.4,border=0.8,arrowsize=2]{->}(0,10)(10,0)
\rput(12,0){$\scriptstyle b$}
\rput(-2,0){$\scriptstyle a$}
\end{pspicture}} \pm\phi(a, b),\]
and keeping the recipe given in Section~\ref{S:SD_hom} for quandle cocycles, one gets braid/knot invariants out of birack/biquandle $2$-cocycles.
\psset{unit=1mm}

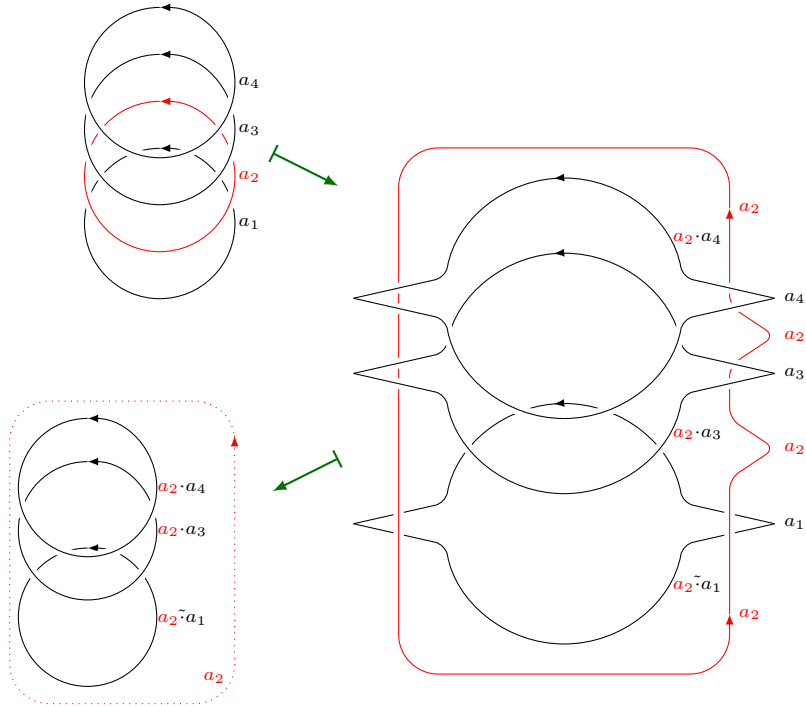
\begin{figure}[h]
\centering
\newsavebox{\UFObox}
\sbox{\UFObox}{
\begin{tikzpicture}[scale=0.23,>=latex]
 \draw [line width=4pt,white] (0,-2.5) circle [radius=4]; 
 \draw [->] (0,1.5) arc [radius=4, start angle=90, end angle= 450];
 \draw [line width=4pt,white] (0,2.5) circle [radius=4]; 
 \draw [->] (0,6.5) arc [radius=4, start angle=90, end angle= 450]; 
 \draw [line width=4pt,white] (0,5) circle [radius=4]; 
 \draw [->] (0,9) arc [radius=4, start angle=90, end angle= 450]; 
 \draw [->,red,rounded corners=15,dotted] (8.5,8) -- (8.5,10) -- (-4.5,10) -- (-4.5,-7.5) -- (8.5,-7.5) -- (8.5,8);
 \node at (3.5,5) [right] {$\scriptstyle{{\color{red}a_2} \cdot a_4}$}; 
 \node at (3.5,2.5) [right] {$\scriptstyle{{\color{red}a_2} \cdot a_3}$}; 
 \node at (8.5,-6) [red,left] {$\scriptstyle{a_2}$}; 
 \node at (3.5,-2.5) [right] {$\scriptstyle{{\color{red}a_2} \wdot a_1}$};   
\end{tikzpicture}
}

\begin{tikzpicture}[scale=0.25,>=latex]
 \draw [line width=4pt,white] (0,-2.5) circle [radius=4]; 
 \draw [->] (0,1.5) arc [radius=4, start angle=90, end angle= 450];
 \draw [line width=4pt,white] (0,0) circle [radius=4]; 
 \draw [->,red] (0,4) arc [radius=4, start angle=90, end angle= 450];
 \draw [line width=4pt,white] (0,2.5) circle [radius=4]; 
 \draw [->] (0,6.5) arc [radius=4, start angle=90, end angle= 450]; 
 \draw [line width=4pt,white] (0,5) circle [radius=4]; 
 \draw [->] (0,9) arc [radius=4, start angle=90, end angle= 450]; 
 \node at (3.7,5) [right] {$\scriptstyle{a_4}$}; 
 \node at (3.7,2.5) [right] {$\scriptstyle{a_3}$}; 
 \node at (3.7,0) [red,right] {$\scriptstyle{a_2}$}; 
 \node at (3.7,-2.5) [right] {$\scriptstyle{a_1}$};    
 \draw [|->, mygreen, thick] (6,1.25)--(9.5,-0.5);
 \draw [|->, mygreen, thick] (9.5,-15)--(6,-16.75);
 \node at (-2,-20) {\usebox{\UFObox}}; 
\end{tikzpicture}
\begin{tikzpicture}[scale=0.4,>=latex]
 \path [draw,line width=4pt,white, rounded corners] (0,1.5) arc (90:170:4) -- (-7,-2.5);
 \path [draw,line width=4pt,white, rounded corners] (0,-6.5) arc (270:190:4) -- (-7,-2.5); 
 \path [draw,line width=4pt,white, rounded corners] (0,-6.5) arc (270:350:4) -- (7,-2.5); 
 \path [draw,line width=4pt,white, rounded corners] (0,1.5) arc (90:10:4) -- (7,-2.5); 
 \draw [->] (0,1.5) arc [radius=4, start angle=90, end angle= 95]; 
 \path [draw, rounded corners] (0,1.5) arc (90:170:4) -- (-7,-2.5);
 \path [draw, rounded corners] (0,-6.5) arc (270:190:4) -- (-7,-2.5); 
 \path [draw, rounded corners] (0,-6.5) arc (270:350:4) -- (7,-2.5); 
 \path [draw, rounded corners] (0,1.5) arc (90:10:4) -- (7,-2.5); 
 \draw [line width=4pt,white,rounded corners=15] (5.5,8) -- (5.5,10) -- (-5.5,10) -- (-5.5,-7.5) -- (5.5,-7.5) -- (5.5,-5.5);
 \draw [line width=4pt,white,rounded corners] (5.5,-5.5) -- (5.5,-1) -- (7,0) -- (5.5,1) -- (5.5,2.75) -- (7,3.75) -- (5.5,4.75) -- (5.5,8);   
 \draw [->,red,rounded corners=15] (5.5,8) -- (5.5,10) -- (-5.5,10) -- (-5.5,-7.5) -- (5.5,-7.5) -- (5.5,-5.5);
 \draw [->,red,rounded corners] (5.5,-5.5) -- (5.5,-1) -- (7,0) -- (5.5,1) -- (5.5,2.75) -- (7,3.75) -- (5.5,4.75) -- (5.5,8);  
 \path [draw,line width=4pt,white, rounded corners] (0,6.5) arc (90:170:4) -- (-7,2.5);
 \path [draw,line width=4pt,white, rounded corners] (0,-1.5) arc (270:190:4) -- (-7,2.5); 
 \path [draw,line width=4pt,white, rounded corners] (0,-1.5) arc (270:350:4) -- (7,2.5); 
 \path [draw,line width=4pt,white, rounded corners] (0,6.5) arc (90:10:4) -- (7,2.5); 
 \draw [->] (0,6.5) arc [radius=4, start angle=90, end angle= 95]; 
 \path [draw, rounded corners] (0,6.5) arc (90:170:4) -- (-7,2.5);
 \path [draw, rounded corners] (0,-1.5) arc (270:190:4) -- (-7,2.5); 
 \path [draw, rounded corners] (0,-1.5) arc (270:350:4) -- (7,2.5); 
 \path [draw, rounded corners] (0,6.5) arc (90:10:4) -- (7,2.5); 
 \path [draw,line width=4pt,white, rounded corners] (0,9) arc (90:170:4) -- (-7,5);
 \path [draw,line width=4pt,white, rounded corners] (0,1) arc (270:190:4) -- (-7,5); 
 \path [draw,line width=4pt,white, rounded corners] (0,1) arc (270:350:4) -- (7,5); 
 \path [draw,line width=4pt,white, rounded corners] (0,9) arc (90:10:4) -- (7,5);    
 \draw [->] (0,9) arc [radius=4, start angle=90, end angle= 95]; 
 \path [draw, rounded corners] (0,9) arc (90:170:4) -- (-7,5);
 \path [draw, rounded corners] (0,1) arc (270:190:4) -- (-7,5); 
 \path [draw, rounded corners] (0,1) arc (270:350:4) -- (7,5); 
 \path [draw, rounded corners] (0,9) arc (90:10:4) -- (7,5);  
 \node at (7,5) [right] {$\scriptstyle{a_4}$}; 
 \node at (3.3,7) [right] {$\scriptstyle{{\color{red}a_2} \cdot a_4}$};  
 \node at (7,2.5) [right] {$\scriptstyle{a_3}$}; 
 \node at (3.3,0.5) [right] {$\scriptstyle{{\color{red}a_2} \cdot a_3}$};   
 \node at (7,0) [red,right] {$\scriptstyle{a_2}$}; 
 \node at (7,3.75) [red,right] {$\scriptstyle{a_2}$}; 
 \node at (5.5,8) [red,right] {$\scriptstyle{a_2}$};  
 \node at (5.5,-5.5) [red,right] {$\scriptstyle{a_2}$}; 
 \node at (7,-2.5) [right] {$\scriptstyle{a_1}$};
 \node at (3.3,-4.5) [right] {$\scriptstyle{{\color{red}a_2} \wdot a_1}$};   
 \node at (3.3,-8.5) {};     
\end{tikzpicture}
\caption{A graphical version of birack cohomology}\label{P:Bir_Cohom}
\end{figure}

Braided and birack cohomology theories were developed in parallel for some time, until it was realized that the guitar map induces an isomorphism between their defining cochain complexes \cite{LebVen}:
\[J^* \colon   (C^{\bullet}_\Bir (S,\ZZ_n),d^{\bullet}_\Bir) \cong (C^{\bullet}_\Br (S,\ZZ_n),d^{\bullet}_\Br).\]
For a biquandle, $J^*$ restricts to an isomorphism between biquandle and braided biquandle cochain complexes. Hence the two cohomology theories are completely identical, and yield the same topological invariants. The proof is once again graphical. Its core is the ``flying saucer'' interpretation of $d^k_\Bir f$. It works as follows. Consider piled circles colored by $a_1, \ldots, a_{k+1}$. Make the $i$th circle inflate or shrink and then disappear, and keep track of the induced color changes; see Fig.~\ref{P:Bir_Cohom} for the inflation situation in the case $n = 4$, $i = 2$. Then evaluate $f$ at the colors of the $k$ remaining circles, and take an alternating sum of these evaluations. 

In our favorite examples, we recover folklore results, namely,
\begin{enumerate}
\item the two forms of group cohomology;
\item the  fact that a rack $(S,\op)$ and its dual $(S,\widetilde{\op})$ have the same cohomology.
\end{enumerate}
For involutive biracks, on the contrary, new results are obtained.

\subsection*{Acknowledgments}
This work was partially supported by the Hamilton Mathematics Institute. The author is grateful to the organizers of the 2017 KIAS Research Station \textit{Self-distributive system and quandle (co)homology theory in algebra and low-dimensional topology} for offering to algebraists and topologists interested in self-distributivity a valuable platform for communication and collaboration.

\bibliographystyle{alpha}
\bibliography{refs}
\end{document}